\documentclass{amsart}
\usepackage{amssymb,euscript,amsmath, mathrsfs, amscd}
\usepackage[dvips]{graphicx}
\usepackage[dvips]{color}

\newcounter{ENUM}
\newcommand{\itm}{\item}
\newenvironment{ilist}{\renewcommand{\theENUM}{\roman{ENUM}}\renewcommand{\itm}{\addtocounter{ENUM}{1}\item[(\theENUM)]}\begin{itemize}\setcounter{ENUM}{0}}{\end{itemize}}
\newenvironment{Ilist}{\renewcommand{\theENUM}{\Roman{ENUM}}\renewcommand{\itm}{\addtocounter{ENUM}{1}\item[(\theENUM)]}\begin{itemize}\setcounter{ENUM}{0}}{\end{itemize}}

\newcommand{\margh}[1]{}

\input xy
\xyoption{all}
\CompileMatrices

\def\cH{{\mathcal H}}

\def\cLH{\mathcal{LH}}

\def\cHom{{\mathcal H}om}
\def\sE{{\mathscr E}}
\def\sF{{\mathscr F}}
\def\sG{{\mathscr G}}

\def\sL{{\mathscr L}}

\def\sO{{\mathscr O}}

\def\vp{\varphi}

\def\id{\operatorname{id}}

\def\im{\operatorname{im}}

\newtheorem{thm}{Theorem}[section]

\newtheorem{lem}[thm]{Lemma}

\theoremstyle{definition}
\newtheorem{defn}[thm]{Definition}

\theoremstyle{remark}
\newtheorem{notn}[thm]{Notation}
\newtheorem{rem}[thm]{Remark}

\numberwithin{equation}{section}
\numberwithin{figure}{section}

\begin{document}
\title{Linked Hom spaces}
\author{Brian Osserman}
\begin{abstract} In this note, we describe a theory of linked Hom spaces
which complements that of linked Grassmannians. Given two chains of vector
bundles linked by maps in both
directions, we give conditions for the space of homomorphisms from one
chain to the other to be itself represented by a vector bundle. We apply
this to present a more transparent version of an earlier construction of 
limit linear series spaces out of linked Grassmannians. 
\end{abstract}

\maketitle

\section{Introduction}

In \cite{os8}, spaces of linked Grassmannians are introduced in order to
use them in a new construction of limit linear series schemes. Given a
base scheme $S$, integers $r<d$, and a chain of vector bundles 
$\sE_1,\dots,\sE_n$ on $S$ of rank $d$, with homomorphisms 
$f_i:\sE_i \to \sE_{i+1}$ and $g_i:\sE_{i+1} \to \sE_i$ satisfying
certain natural conditions, the associated linked Grassmannian
parametrizes tuples of subbundles $\sF_i \subseteq \sE_i$ of rank $r$,
which are all mapped into one another under the $f_i$ and $g_i$. These
schemes behave like flat degenerations of the classical Grassmannian 
$G(r,d)$, and indeed according to \cite{o-h1}, whenever the $f_i$ and $g_i$ 
are generically isomorphisms, the linked Grassmannian does in fact yield
a flat degeneration of $G(r,d)$.

In this note, we study a Hom version of the linked Grassmannian 
construction, which arises naturally in the construction of limit linear
series spaces out of linked Grassmannians. Given chains of vector bundles
$\sF_1,\dots,\sF_n$ and $\sG_1,\dots,\sG_n$ on $S$, with the $\sF_i$ of rank 
$r$ and the $\sG_i$ of rank $m$, and homomorphisms
$f_i:\sF_i \to \sF_{i+1}$, $f^i:\sF_{i+1} \to \sF_i$,
$g_i:\sG_i \to \sG_{i+1}$, $g^i:\sG_{i+1} \to \sG_i$, we determine
natural conditions (see Definition \ref{defn:link-hom}
below) so that the resulting linked Hom functor 
$\cLH(\{\sF_i\},\{\sG_i\},\{f_i,f^i\},\{g_i,g^i\})$ parametrizing
tuples of morphisms $\vp_i:\sF_i \to \sG_i$ which commute with all
$f_i,f^i,g_i,g^i$ is well behaved. Specifically, we show:

\begin{thm}\label{thm:main} A linked Hom functor 
$\cLH(\{\sF_i\},\{\sG_i\},\{f_i,f^i\},\{g_i,g^i\})$ 
is represented by a vector bundle $LH$ on $S$ of rank $rm$.
\end{thm}

Using this result, we are able to give a less ad hoc and more symmetric
version of the construction of limit linear series spaces out of linked 
Grassmannians. Although the resulting shift in perspective is relatively 
minor, we believe it will be important for a new application to higher-rank
Brill-Noether theory. Specifically, Bertram, Feinberg and Mukai observed
that in the case of rank $2$ vector bundles of canonical determinant, there
are symmetries which result in the spaces of bundles with prescribed numbers
of sections having dimension higher than previously expected, when they
are non-empty. One would like to show such spaces (and generalizations,
as suggested in \cite{os16}) are nonempty via degeneration techniques, but 
the limit linear series construction introduced by Eisenbud and Harris 
in \cite{e-h1} and generalized to higher rank by Teixidor i Bigas in
\cite{te1} does not appear 
superficially to preserve the necessary symmetries introduced by the 
canonical determinant condition. The simpler nature of our construction 
should lead to a proof that the dimensions of spaces of limit linear
series with special determinants satisfy the desired modified lower
bounds.

\subsection*{Acknowledgements} I would like to thank Montserrat Teixidor
i Bigas for helpful comments.

\section{Linked Hom spaces}

We work throughout over a fixed base scheme $S$.

The basic definition is as follows:

\begin{defn}\label{defn:basic} Let $r,m,n$ be integers. 
Suppose that $\sF_1,\dots,\sF_n$ are vector bundles of rank $r$ on $S$
and $\sG_1,\dots,\sG_n$ are vector bundles of rank $m$ on $S$. Suppose
also that we have homomorphisms 
$$f_i:\sF_i \to \sF_{i+1},\quad
f^i:\sF_{i+1} \to \sF_i,\quad g_i:\sG_i \to \sG_{i+1},\quad
g^i:\sG_{i+1} \to \sG_i$$
for each $i=1,\dots,n-1$. Then the functor
$$\cLH(\{\sF_i\},\{\sG_i\},\{f_i,f^i\},\{g_i,g^i\})$$
associates to an
$S$-scheme $T$ the tuples of homomorphism $\vp_i:\sF_i|_T \to \sG_i|_T$ such
that $\vp_{i+1} \circ f_i =g_i \circ \vp_i$ and 
$\vp_i \circ f^i = g^i \circ \vp_{i+1}$ for $i=1,\dots,n-1$.
\end{defn}

\begin{lem}\label{lem:representable} In the situation of Definition
\ref{defn:basic}, the functor 
$$\cLH(\{\sF_i\},\{\sG_i\},\{f_i,f^i\},\{g_i,g^i\})$$
is represented by
an affine scheme $LH$ of finite presentation over $S$.
\end{lem}

Note that the functor naturally takes values in $\sO_T$-modules, and as
a result $LH$ has the structure of a (not necessarily quasi-coherent) 
$\sO$-module over the big Zariski (or etale, or fppf) site on $S$. 

\begin{proof} It is clear that our functor is a subfunctor of 
$\bigoplus_{i=1}^n \cHom(\sF_i,\sG_i)$, which is represented by a vector
bundle over $S$. Moreover, it is cut out by the conditions
$\vp_{j+1} \circ f_j -g_j \circ \vp_j=0$ and 
$\vp_j \circ f^j - g^j \circ \vp_{j+1}=0$, which can be viewed as
the preimage of the zero section under morphisms
$$\bigoplus_{i=1}^n \cHom(\sF_i,\sG_i) \to \cHom(\sF_j,\sG_{j+1})$$
and
$$\bigoplus_{i=1}^n \cHom(\sF_i,\sG_i) \to \cHom(\sF_{j+1},\sG_j)$$
respectively, and thus give closed (and finitely generated) conditions. 
We thus conclude that $LH$ exists, with the desired properties.
\end{proof}

We now describe the additional conditions we will impose in order to
obtain good behavior for $LH$. Rather than literally considering a Hom 
version of linked Grassmannians (that is, imposing the same conditions on
both the $\sF_i$ and the $\sG_i$ as were imposed on the $\sE_i$ in
a linked Grassmannian), we will impose somewhat stronger conditions on
the $\sG_i$ and weaker ones on the $\sF_i$. This not only leads to good
behavior, but is precisely the situation which arises in the application
to construction of limit linear series spaces. See Remark \ref{rem:compare}
below for further discussion.

\begin{defn}\label{defn:link-hom} In the situation of Definition 
\ref{defn:basic} and Lemma
\ref{lem:representable}, we say that $LH$ is a {\bf linked Hom space}
if the following conditions are satisfied:
\begin{Ilist}
\itm There exists some $s \in \Gamma(S,\sO_S)$ such that 
$$f_i f^i=f^i f_i=g_i g^i = g^i g_i=s \id$$
for all $i=1,\dots,n-1$.
\itm For all $x \in S$ with $s=0$ in $\kappa(x)$, and all $i=1,\dots,n-1$, 
the kernel of $g_i$ 
restricted to $\sG_i \otimes \kappa(x)$ is precisely the
image of $g^i$, and vice versa.
\itm For all $x \in S$ with $s=0$ in $\kappa(x)$, and all $i=1,\dots,n-2$,
we have $\im g_i$ complementary to $\ker g_{i+1}$, and $\im g^{i+1}$ 
complementary to $\ker g^i$ in $\sG_{i+1} \otimes \kappa(x)$.
\end{Ilist}
\end{defn}

The following notation will be convenient:

\begin{notn} In the situation of Definition \ref{defn:link-hom},
given $i<j$, denote by $f_{i,j}$, $f^{j,i}$, $g_{i,j}$ and
$g^{j,i}$ the compositions
$$f_{j-1} \circ f_{j-2} \circ \cdots \circ f_{i+1} \circ f_i,$$ 
$$f^i \circ f^{i+1} \circ \cdots \circ f^{j-2} \circ f^{j-1},$$
$$g_{j-1} \circ g_{j-2} \circ \cdots \circ g_{i+1} \circ g_i,$$ 
and
$$g^i \circ g^{i+1} \circ \cdots 
\circ g^{j-2} \circ g^{j-1},$$
respectively. We use the convention that 
$f_{i,i}$, $f^{i,i}$, $g_{i,i}$ and $g^{i,i}$ all denote the identity map.
\end{notn}

Our hypotheses on the chain $\sG_i$ yield the following simple structure:

\begin{lem}\label{lem:structure} In the situation of Definition
\ref{defn:basic}, suppose that $LH$ is a linked Hom space. Then locally
on $S$, we have decompositions 
$$\sG_i \cong \sG_i' \oplus \sG_i''$$
with $\sG_i'$ and $\sG_i''$ free $\sO_S$-modules,
satisfying the following conditions for all $i=1,\dots,n-1$:
\begin{ilist}
\itm $g_i$ maps $\sG_i'$ to $\sG_{i+1}'$ and $\sG_i''$ to $\sG_{i+1}''$,
and similarly for $g^i$.
\itm The induced maps $(g_i)':\sG'_i \to \sG'_{i+1}$ and 
$(g^i)'':\sG''_{i+1} \to \sG''_i$ are isomorphisms.
\end{ilist}
\end{lem} 

\begin{proof} Given $x \in S$, if $s \neq 0$ in $\kappa(x)$, then there
exists an open neighborhood $U$ of $x$ on which $s$ is a unit, and by 
condition (I) of a linked Hom space, the $g_i$ and $g^i$ are isomorphisms. 
We may thus set $\sG_i'=\sG_i$ and $\sG_i''=0$ for all $i$.

On the other hand, if $s=0$ in $\kappa(x)$, choose subspaces 
$\bar{\sG}'_1 \subseteq \sG_1 \otimes \kappa(x)$ and 
$\bar{\sG}''_n \subseteq \sG_n \otimes \kappa(x)$,
complementary to $\ker g_1$ and $\ker g^{n-1}$ respectively. Putting
together conditions (II) and (III) for a linked Hom space, we see
that for all $i$, we have $g_{1,i}$ injective on $\bar{\sG}'_1$,
and $g^{n,i}$ injective on $\bar{\sG}''_n$, and furthermore
setting $\bar{\sG}'_i = g_{1,i} \bar{\sG}'_1$ and
$\bar{\sG}''_i = g^{n,i} \bar{\sG}''_n$ we have that
$\bar{\sG}'_i$ and $\bar{\sG}''_i$ are complementary in
$\sG_i \otimes \kappa(x)$. We thus obtain a decomposition of
the desired form over $\kappa(x)$. In particular, if 
$\dim \bar{\sG}'_1 = \ell$, then $\dim \bar{\sG}''_n = m-\ell$.

Now, choose $U$ any affine open neighborhood of $x$ on which the 
$\sG_i$ become free. Multiplying by units in $\sO_{S,x}$ as necessary,
we may choose bases for $\bar{\sG}'_1$ and $\bar{\sG}''_n$ which can
be lifted to some $v_1,\dots,v_{\ell} \in \sG_1(U)$ and 
$w1_,\dots,w_{m-\ell} \in \sG_n(U)$.
Let $\sG'_1$ and $\sG''_n$ be the submodules of $\sG_1|_U$ and 
$\sG_n|_U$ spanned by the $v_i$ and $w_i$, respectively.
Restrict $U$ as necessary to the complement of the closed subschemes on
which the maps
$$\sO_U^{\oplus \ell} \overset{v_i}{\to} \sG_1$$
$$\sO_U^{\oplus m-\ell} \overset{w_i}{\to} \sG_n$$
and for each $i=1,\dots,n$
$$g_{1,i} \oplus g^{n,i}: \sG'_1 \oplus \sG''_n \to \sG_i$$
do not have full rank.
To see that the last condition make sense, we observe that under the
hypothesis that the first map has full rank, it must in fact
be an isomorphism onto its image $\sG'_1$, and similarly for $\sG''_n$. 
It suffices to check this on stalks, and the full rank condition implies 
that on each stalk, the fiber of $\sG'_1$ has dimension $\ell$. Thus
Nakayama's lemma implies that $\sG'_1$ is free of rank $\ell$ on stalks,
and thus the map induced by the $v_i$ is an isomorphism on stalks, as 
desired. The same argument goes through for $\sG''_n$.
 
As before, we set $\sG'_i = g_{1,i} \sG'_1$ and
$\sG''_i = g^{n,i} \sG''_n$, and we claim that we obtain the desired
decomposition. Our first claim is that $g_{1,i}:\sG'_1 \to \sG'_i$ is
an isomorphism for each $i$, and similarly for $g^{n,i}$. In particular,
$\sG'_i$ and $\sG''_i$ are free $\sO_U$-modules, and condition (ii) is
satisfied. Arguing as before,
this follows from Nakayama's lemma and the full rank hypotheses.
Next, for each $i$, we have the natural morphism
$$\sG_i' \oplus \sG_i'' \to \sG_i$$
induced by the inclusions, which we wish to show is an isomorphism.
Since both modules are free of rank
$m$, to obtain the desired decomposition it suffices to check surjectivity, 
which holds on fibers by the full rank hypothesis, and hence on stalks
by Nakayama's lemma.

Since we have already checked condition (ii), it remains to check that (i)
is satisfied. We see that
$g_i$ maps $\sG'_i$ to $\sG'_{i+1}$ by construction. On the other
hand, $g_i \sG''_i = g_i g^i \sG''_{i+1} = s \sG''_{i+1}$ by condition (I)
of a linked Hom space. Thus, $g_i$
preserves the decomposition, and we see
similarly that the maps $g^i$ preserve the decomposition, giving us
condition (i), as desired.
\end{proof}

We can now prove our main theorem:

\begin{proof}[Proof of Theorem \ref{thm:main}] The question being local
on $S$, we may assume without loss of generality that the $\sG_i$ 
decompose as in Lemma \ref{lem:structure}. Denote by $m_1$ the rank of
$\sG'_1$ and by $m_2$ the rank of $\sG''_n$. The existence of the isomorphisms
$(g_i)'$ and $(g^i)''$ imply that $\sG'_i$ has rank $m_1$ for all 
$i$, and $\sG''_i$ has rank $m_2$ for all $i$. In particular, $m_1+m_2=m$. 
We then claim that in fact $LH$ represents 
$$\cH:= \cHom(\sF_1,\sG''_1) \oplus \cHom(\sF_n,\sG'_n)$$ 
and is thus a vector bundle of rank $rm_1+rm_2=rm$, as desired.

For $i<j$, denote by $(g_{i,j})'$ and
$(g^{j,i})''$ the compositions
$$(g_{j-1})' \circ (g_{j-2})' \circ \cdots \circ (g_{i+1})' \circ (g_i)'$$ 
and
$$(g^i)'' \circ (g^{i+1})'' \circ \cdots \circ (g^{j-2})'' \circ (g^{j-1})'',$$
respectively. We use the convention that $(g_{i,i})'$ and 
$(g^{i,i})''$ denote the identity map.

Clearly, there is a natural forgetful map from $h_{LH}$ to $\cH$. 
To construct the inverse map, if we are given 
$\vp_1'' \in \cHom(\sF_1,\sG''_1)_T$
and 
$\vp_n' \in \cHom(\sF_n,\sG'_n)_T$,
we can construct $\vp_i \in \cHom(\sF_i,\sG_i)_T$ for all $i$ by
setting 
$$\vp_i=((g_{i,n})')^{-1} \circ \vp_n' \circ f_{i,n} \oplus
((g^{i,1})'')^{-1} \circ \vp_1'' \circ f^{i,1}.$$
We check the linkage condition directly: we have
\begin{align*} \vp_{i+1} \circ f_i & 
=((g_{{i+1},n})')^{-1} \circ \vp_n' \circ f_{i+1,n} \circ f_i 
\oplus ((g^{i+1,1})'')^{-1} \circ \vp_1'' \circ f^{i+1,1} \circ f_i \\
& =((g_{{i+1},n})')^{-1} \circ \vp_n' \circ f_{i,n}
\oplus s ((g^{i+1,1})'')^{-1} \circ \vp_1'' \circ f^{i,1},\end{align*} 
while
\begin{align*} g_i \circ \vp_i &
= g_i \circ ((g_{i,n})')^{-1} \circ \vp_n' \circ f_{i,n} 
\oplus g_i \circ ((g^{i,1})'')^{-1} \circ \vp_1'' \circ f^{i,1} \\
& = ((g_{i+1,n})')^{-1} \circ \vp_n' \circ f_{i,n} 
\oplus g_i \circ (g^i)'' \circ 
((g^{i+1,1})'')^{-1} \circ \vp_1'' \circ f^{i,1},\end{align*}
so we conclude equality from the fact that $g_i \circ (g^i)'' = s \id$.
Similarly, we have
$\vp_i \circ f^i = g^i \circ \vp_{i+1}$, so the linkage condition holds.
Moreover, the fact that $(g_i)'$ and $(g^i)''$ are isomorphisms 
together with the linkage condition imply that
the $\vp_i$ are uniquely determined by $\vp_1''$ and $\vp_n'$,
so we obtain the desired isomorphism of functors and conclude the theorem.
\end{proof}

\begin{rem}\label{rem:compare} 
Note that the conditions on the $f_i,f^i$ are quite minimal, and the
only additional conditions we impose on the $g_i,g^i$ beyond that for
a linked Grassmannian are that $\im g_i$ actually be complementary
to $\ker g_{i+1}$ and similarly for the $g^i$, while for linked 
Grassmannians we required only that $\im g_i \cap \ker g_{i+1}=(0)$
and $\im g^{i+1} \cap \ker g^i=(0)$. 

We see moreover that this mild strengthening is indeed necessary in order 
for Theorem \ref{thm:main} to hold: consider the case that $r=1$, $m=3$, 
and $n=3$, with $s=t^2$, all modules free, and in terms of chosen bases,
maps given as follows: 
$$f_1=f_2=f^1=f^2=t,$$
$$g_1=\begin{bmatrix} 1 & 0 & 0 \\ 0 & s & 0 \\ 0 & 0 & s \end{bmatrix},
\quad
g^1=\begin{bmatrix} s & 0 & 0 \\ 0 & 1 & 0 \\ 0 & 0 & 1 \end{bmatrix},$$
$$g_2=\begin{bmatrix} 1 & 0 & 0 \\ 0 & 1 & 0 \\ 0 & 0 & s \end{bmatrix},
\quad
g^2=\begin{bmatrix} s & 0 & 0 \\ 0 & s & 0 \\ 0 & 0 & 1 \end{bmatrix}.$$
First suppose that $S$ is a point, and that $t=0$. In this case, the space 
$LH$ parametrizes triples of 
vectors in $\sG_1,\sG_2,\sG_3$ which are in the kernels of all relevant
maps, and this has dimension $4$, which is larger than $rm=3$.

If instead $S$ is the spectrum of a DVR, and $t$ is a uniformizer, then
we have that $LH$ has dimension $3$ over the generic point, but 
dimension $4$ over the closed point, so $LH$ is not even flat over $S$.
\end{rem}

\section{Application to limit linear series}

We conclude with a brief explanation of the application of Theorem
\ref{thm:main} to the construction of spaces of limit linear series.
We follow the notation of \cite{os8}. We begin with a family 
$\pi:X \to B$ of curves of genus $g$
satisfying the conditions of a smoothing family, as
specified in Definition 3.1 of {\it loc. cit.} In particular,
over the locus $\Delta$ of $B$ over which $X$ is not smooth, $X$ is
nodal, consisting of two smooth components $Y$ and $Z$ glued along a
section $\Delta'$ of $\pi$. Although \cite{os8} also treats the case in which
$\pi$ is smooth, we will only be interested in the other two cases,
where $\Delta=B$ or where $\Delta$ is a Cartier divisor in $B$. 
Given integers $r,d$, Definition 4.5 of {\it loc. cit.} then describes a 
limit linear series functor which can be summarized roughly as follows.

In the case that $\Delta=B$, to a $B$-scheme $T$, the functor associates 
the set of line bundles $\sL$ on $X|_T$ of degree $d$ on $Y|_T$ and degree
$0$ on $Z|_T$, together with rank-$(r+1)$ subbundles $V_i$ of $\pi_* \sL_i$,
where $\sL_i$ is obtained from $\sL$ by gluing together the line bundles
$\sL|_Y (-i\Delta')$ and $\sL|_Z (i\Delta')$. We further require that each
$V_i$ map into $V_{i+1}$ under the maps which are zero on $Y$ and the
natural inclusion on $Z$, and vice versa.

In the case that $\Delta$ is a Cartier divisor in $B$, to a $B$-scheme $T$
the functor associates the set of line bundles $\sL$ on $X|_T$ of degree $d$,
having degree $d$ on $Y|_T$ and degree $0$ on $Z|_T$, together with 
rank-$(r+1)$ subbundles $V_i$ of $\pi_* \sL_i$, where $\sL_i:=\sL \otimes
\sO_X(iY)|_T$.
We further require that the $V_i$ map into one another under the maps which 
are the natural inclusions in one direction, and which are obtained from
a choice of isomorphism $\sO_X(-Z) \cong \sO_X(Y)$ in the other.

In fact, the definition also allows for imposed ramification along 
sections, but this part of the construction is unaffected by Theorem
\ref{thm:main}, so for the sake of simplicity we ignore it. The 
construction of the scheme representing the functor (as given in
Theorem 5.3 of \cite{os8}) then proceeds as
follows. Let $P$ be the relative Picard scheme parametrizing line bundles
of degree $d$ on $Y$ and degree $0$ on $Z$, and let $\sL$ be the
universal line bundle on $X \times_B P$. Let $\sL_i$ be obtained from
$\sL$ as in the definition of the functor, depending on which case we
are in. Then the $\sL_i$ also have maps in both directions as described 
above. Choose $D$ a sufficiently ample divisor on $X$, and consider
the linked Grassmannian $LG$ parametrizing tuples of rank-$(r+1)$ subbundles
$V_i$ of $p_{2*} \sL_i(p_1^* D)$ which map into one another under the given
maps. This space has dimension $\dim B + g +(r+1)(d+\deg D-g-r)$. We can cut
out the scheme representing our functor by requiring that the $V_i$
are in fact contained in $p_{2*} \sL_i$ for each $i$, or equivalently that
the induced maps $V_i \to p_{2*} (\sL_i(p_1^* D)|_{p_1^* D})$ are all zero. 
However, the trick is to give the appropriate bounds on the dimension. In 
\cite{os8}, this was accomplished in an {\it ad hoc} manner by writing 
$D=D_Y+D_Z$ for divisors $D_Y,D_Z$ supported entirely on $Y$ and $Z$ 
respectively, and checking that imposing vanishing along $p_1^* D$ for all 
$i$ is equivalent to imposing vanishing on $p_1^* D_Y$ for $V_0$, and 
vanishing on $p_1^* D_Z$ for $V_d$. 

With Theorem \ref{thm:main}, we can proceed more directly. We can rephrase
the above discussion as follows: our space $LG$ carries the universal 
chains of bundles
$\{V_i\}$ and $\{p_{2*} (\sL_i(p_1^* D)|_{p_1^* D})\}$, of ranks $r+1$
and $\deg D$ respectively, together with a linked homomorphism between 
them, and the desired limit linear series space is cut out precisely by 
the condition that this linked homomorphism be $0$ for all $i$. Now, it 
is easy to see that our chains of bundles
satisfy the conditions for a linked Hom space. Indeed, condition (I)
is inherited from the $p_{2*} \sL_i(p_1^* D)$. Next, we have 
that the maps between 
the $p_{2*} (\sL_i(p_1^* D)|_{p_1^* D})$ fail to be isomorphisms 
precisely over $\Delta$. Over $\Delta$, we have that a section of 
$\sL_i(p_1^* D)|_{p_1^* D}$
is in the image of $\sL_{i-1}(p_1^* D)|_{p_1^* D}$ if and only if 
it vanishes along $p_1^*(D \cap Y)$, which is to say if it is in the
kernel of the map back to $\sL_{i-1}(p_1^* D)|_{p_1^*D}$. Thus, condition
(II) is satisfied. On the other hand, a section is in the kernel of the map 
to $\sL_{i+1}(p_1^* D)|_{p_1^* D}$ if and only if it vanishes along 
$p_1^*(D \cap Z)$. However, $\sL_i(p_1^* D)|_{p_1^* D}$ decomposes
as a direct sum of the spaces of sections vanishing on $p_1^*(D \cap Y)$
and on $p_1^*(D \cap Z)$, so we obtain condition (III) as well. Because
the limit linear series space is the preimage of the zero section under 
the given map to the linked Hom space, we conclude from Theorem 
\ref{thm:main} that each component has codimension at most $(r+1) \deg D$ 
inside $LG$, and thus that each component has dimension at least
$$\dim B+ g + (r+1)(d-g-r)=\dim B + (r+1)(d-r) - rg,$$
agreeing with the formula obtained in Theorem 5.3 of \cite{os8}.

\bibliographystyle{hamsplain}
\bibliography{hgen}

\end{document}